\documentclass[12pt]{article}
\usepackage{amsmath,amssymb,amsthm,xspace,picinpar,graphicx,textcomp}
\newcommand{\old}[1]{}

\renewcommand{\th}{th\xspace}

\newcommand{\textem}[1]{\textit{#1\/}}
\newtheorem{thm}{Theorem}[subsubsection]
\newtheorem{theorem}{Theorem}[subsubsection]

\newtheorem{lemma}[theorem]{Lemma}

\newtheorem{prop}[theorem]{Proposition}

\newcommand{\R}{{\mathbb R}}

\newcommand{\E}{{\mathbb E}}

\newcommand{\eps}{\varepsilon}

\newcommand{\Var}{{\rm Var}}

\newcommand{\vs}{\vspace{.15 in}}

\begin{document}
\title{Random-Turn Hex and Other Selection Games}
\old{
\author{Yuval Peres\thanks{U.C. Berkeley.  Research supported in part by NSF grants \#DMS-0244479 and \#DMS-0104073.}\ \thanks{Microsoft Research.}
\and Oded Schramm\footnotemark[2]
\and Scott Sheffield\footnotemark[1] \footnotemark[2]
\and David B. Wilson\footnotemark[2]}
}
\author{Yuval Peres, Oded Schramm, Scott Sheffield, and David B. Wilson}
\date{} \maketitle

\old{
\begin{abstract}
  The game of Hex has two players who take turns placing
  stones of their respective colors on the hexagons
  of a rhombus-shaped hexagonal
  grid.  Black wins by completing a crossing between two opposite
  edges, while White wins by completing a crossing between the other
  pair of opposite edges.  Although ordinary Hex is famously difficult
  to analyze, \textit{Random-Turn\/} Hex---in which players toss a coin
  before each turn to decide who gets to place the next stone---has a
  simple optimal strategy.  It belongs to a general class of
  random-turn games---called selection games---in which the
  expected payoff when both players play the random-turn game
  optimally is the same as when both players play randomly.  We also
  describe the optimal strategy and study the expected length of the
  game under optimal play for Random-Turn Hex and several other
  selection games.
\end{abstract}
}

\renewcommand{\thesubsubsection}{\arabic{subsubsection}}
\let\section\subsubsection
\let\subsection\paragraph
\section{INTRODUCTION.}
\subsection*{Overview.}

\newcommand{\olympcapt}{A game of Hex played at the 5\th Computer
  Olympiad in London, August 24, 2000 \cite{hexinfo}.  Queenbee
  \cite{queenbee} (black) played against Hexy \cite{hexy} (white).
  (There is a special rule regarding the first move whose purpose is
  to offset the advantage of the first player, see~\cite{hexinfo}.)
 \label{fig:olympiad}}
The game of \textem{Hex}, invented independently by Piet Hein in 1942
and John Nash in 1948 \cite{nash}, has two players who take turns placing stones
of their respective colors on the hexagons of a rhombus-shaped hexagonal grid
(see Figure~\ref{fig:olympiad}).
A player wins by completing a path connecting the two opposite sides
of his or her color.
Although it is easy to show that player~I has a winning strategy,
it is far from obvious what that strategy is.  Hex is usually played
on an $11 \times 11$ board (e.g., the commercial version by Parker
Brothers$^{\text{\tiny\textregistered}}$ and the Computer Olympiad Hex Tournament both use $11\times
11$ boards),
for which the optimal strategy is not yet known.  (For a book on practical Hex strategy see \cite{browne}, and for further information on Hex see \cite{stewart}, \cite{gardner}, \cite{hexinfo}, or \cite{gale}.)
\begin{figure}[phtb]
\centerline{\includegraphics[width=0.5\textwidth]{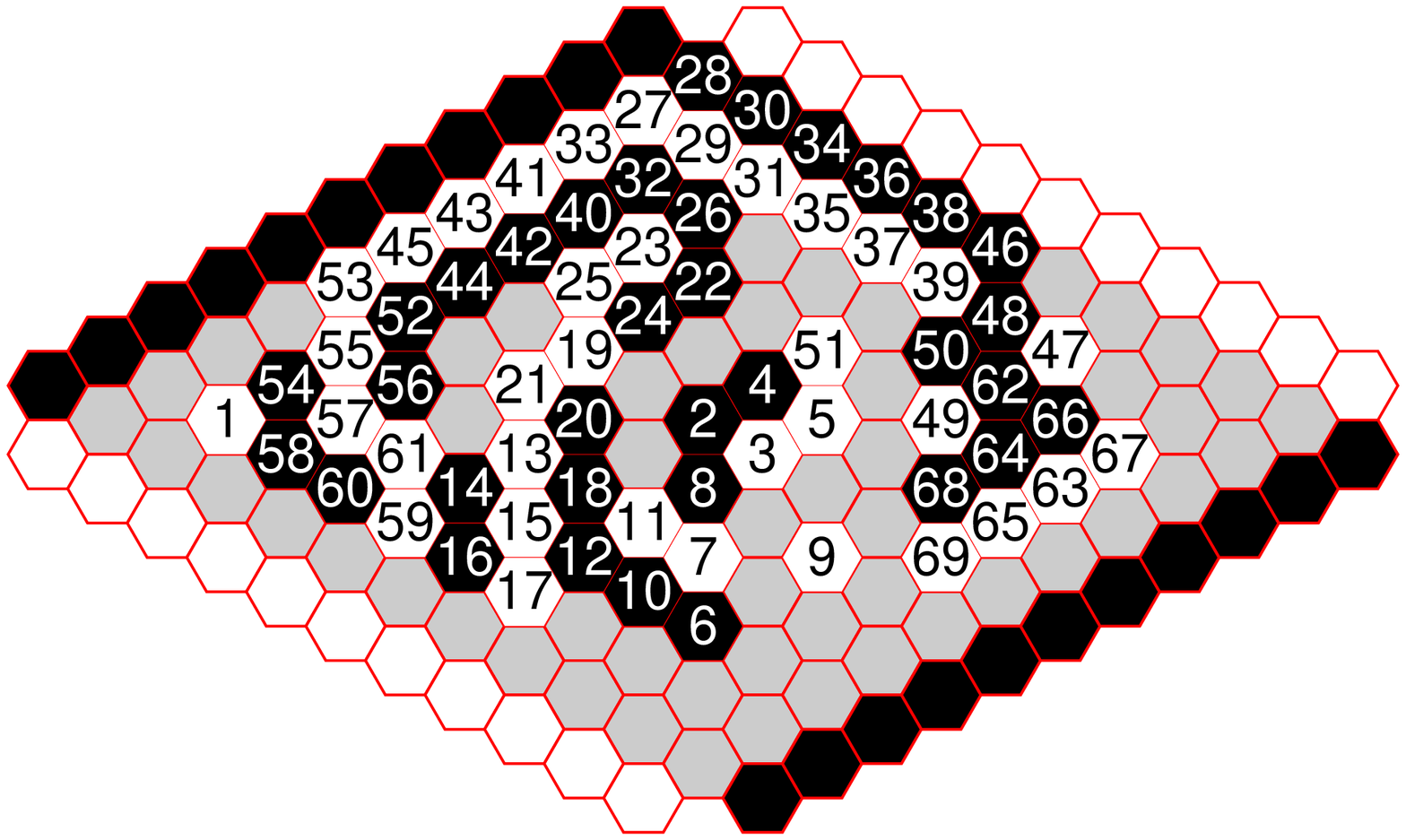}}
\caption{\olympcapt}
\end{figure}

\textem{Random-Turn Hex} is the same as ordinary Hex, except that
instead of alternating turns, players toss a coin before each turn
to decide who gets to place the next stone.
Although ordinary Hex is famously difficult to analyze, the optimal
strategy for Random-Turn Hex turns out to be very simple.
We introduce random-turn games in part because they are in many cases
more tractable than their deterministic analogs; they also exhibit
surprising structure and symmetry.  In one general class of games
called \textit{selection games}, which includes Hex, the probability
that player~I wins when both players play \textit{optimally\/} is
precisely the probability that player~I wins when both players
play \textit{randomly}.  Combining this with Smirnov's recent result
about percolation on the hexagonal lattice \cite{smirnov},
we will see that in a
certain ``fine lattice limit'' winning probabilities of Random-Turn
Hex are a ``conformal invariant'' of the board shape.  In this and
other games, the set of moves played during an entire game (when both
players play optimally) has an intriguing fractal structure.

On a philosophical note, random-turn games are natural models for
real-world conflicts, where opposing agents (political parties,
lobbyists, businesses, militaries, etc.)\ do not alternate turns.
Instead, they continually seek to improve their positions incrementally.
The extent to which they succeed is (at least partially)
unpredictable and may be modeled using randomness.  Random-turn games
are also closely related to ``Richman games'' \cite{MR1685133}.  In a
Richman game, each player offers money to the other player for the
right to make the next move, and the player who offers more gets to
move.  (At the end of the Richman game, the money has no value.)
Another class of random-turn games, where the moves of one player may
be reversed by the other player, is studied in
\cite{MR1685133} and \cite{pssw-tug}.

For simplicity, we limit the discussion in this article to two-player,
zero-sum, random-turn games with perfect information.

\subsection*{Random-turn selection games.}
Now we describe a general class of games that includes the famous game of Hex.
Let $S$ be an $n$-element set, which will sometimes be called the
\textem{board}, and let $f$ be a function from the $2^n$ subsets of $S$ to
$\mathbb R$.  A \textem{selection game} is played as follows: the
first player selects an element of $S$, the second player selects one
of the remaining $n-1$ elements, the first player selects one of the
remaining $n-2$, and so forth, until all elements have been chosen.
Let $S_1$ and $S_2$ signify the sets chosen by the first and second
players, respectively.  Then player~I receives a payoff of $f(S_1)$
and player~II a payoff of $ - f(S_1)$ (selection games are zero-sum).
The following are examples of selection games:

\newcommand{\spar}[1]{\subsection*{#1.}}

\spar{Hex}
Here $S$ is the set of hexagons on a rhombus-shaped $L
\times L$ hexagonal grid, and $f(S_1)$ is $1$ if $S_1$ contains a
left-right crossing, $-1$ otherwise (see Figures~\ref{fig:olympiad}
and \ref{fig:hex11}).  In this case, once $S_1$ contains a left-right
crossing or $S_2$ contains an up-down crossing (which precludes the
possibility of $S_1$ having a left-right crossing), the outcome is
determined and there is no need to continue the game.
\medskip

We also consider Hex played on other types of boards.
In the general setting, some hexagons are given to the first or second
players before the game has begun.  One of the reasons for considering
such games is that after a number of moves are played in ordinary Hex,
the remaining game has this form.

\spar{Surround}
The famous game of \textit{Go\/} is not a selection game (for one, a player
can remove an opponent's pieces), but the game of \textit{Surround}, in
which, as in Go, surrounding area is important, is a selection game.
In this game $S$ is the set of $n$ hexagons in a hexagonal grid
(of any shape).  At the end of the game, each hexagon is recolored to
be the color of the outermost cluster surrounding it (if there is such
a cluster).  The payoff $f(S_1)$ is the number of hexagons recolored
black minus the number of hexagons recolored white.  (Another natural
payoff function is $f^*(S_1)=\operatorname{sign}(f(S_1))$.)

\spar{Full-Board Tic-Tac-Toe}
Here $S$ is the set of spaces in a $3
\times 3$ grid, and $f(S_1)$ is the number of horizontal, vertical, or
diagonal rows in $S_1$ minus the number of horizontal, vertical, or
diagonal rows in $S \backslash S_1$.  This is different from ordinary
tic-tac-toe in that the game does not end after the first row is
completed.

\spar{Team Captains}
Two team captains are choosing
baseball teams from a finite set $S$ of $n$ players for the purpose of
playing a single game against each other.  The payoff $f(S_1)$ for the
first captain is the probability that the players in $S_1$ (together
with the first captain) would beat the players in $S_2$ (together with
the second captain).  The payoff function may be very complicated
(depending on which players are skilled at which positions, which players have
played together before, which players get along well with which
captain, etc.).  Because we have not specified the payoff function,
this game is as general as the class of selection games.

\vs

Every selection game has a random-turn variant in which at each turn
a fair coin is tossed to decide who moves next.

\vs

Consider the following questions:

\begin{enumerate}
\item What can one say about the probability distribution of $S_1$
  after a typical game of optimally played Random-Turn Surround?
\item More generally, in a generic random-turn selection game, how
  does the probability distribution of the final state depend on the
  payoff function $f$?
\item Less precise: Are the teams chosen in Random-Turn Team Captains
  ``good teams'' in any objective sense?
\end{enumerate}

\noindent
The answers are surprisingly simple.

\section{OPTIMAL STRATEGY.}
\label{c.opt}

A (pure) \textem{strategy} for a given player
in a random-turn selection game is a
map $M$ from pairs of disjoint subsets $(T_1, T_2)$ of $S$ to elements
of $S$.  Here $M(T_1, T_2)$ indicates the element that the player will
pick if given a turn at a time in the game at which player~I
has thus far picked the elements of $T_1$ and player~II has picked
the elements of $T_2$.

Denote by $E(T_1, T_2)$ the expected payoff for player~I at this
stage in the game, assuming that both players play optimally with the
goal of maximizing expected payoff.  As is true for all finite
perfect-information, two-player games, $E$ is well defined, and one can
compute $E$ and the set of possible optimal
strategies inductively as follows.  First, if $T_1 \cup T_2 = S$, then
$E(T_1, T_2) = f(T_1)$.  Next, suppose that we have computed $E(T_1, T_2)$
whenever $|S \backslash (T_1 \cup T_2)| \leq k$.  Then if $|S
\backslash (T_1 \cup T_2)|=k+1$, and player~I has the chance to
move, player~I will play optimally if and only if she chooses an $s$
from $S \backslash (T_1 \cup T_2)$ for which $E(T_1 \cup \{s\}, T_2)$ is
maximal.  (If she chose any other $s$, this would reduce her expected
payoff.)  Similarly, player~II plays optimally if and only if she
minimizes $E(T_1, T_2 \cup \{s\})$ at each stage.

The foregoing analysis also demonstrates a well-known fundamental fact
about finite, turn-based, perfect-information games: both players have
optimal pure strategies (i.e., strategies that do not require flipping
coins), and knowing the other player's strategy does not give a player
any advantage when both players play optimally.  (This contrasts with
the situation in which the players play ``simultaneously,'' as they do
in Rock-Paper-Scissors.)

\begin{theorem} \label{t.rts}
  The value of a random-turn selection game is the expectation of
  $f(T)$ when a set $T$ is selected randomly and uniformly among
  all subsets of $S$.  Moreover, any optimal strategy for one of the
  players is also an optimal strategy for the other player.
\end{theorem}

\begin{proof}
  If player~II plays any optimal strategy, player~I can achieve the
  expected payoff $\E[f(T)]$ by playing exactly the same strategy
  (since, when both players play the same strategy, each element will
  belong to $S_1$ with probability $1/2$, independently).  Thus, the
  value of the game is at least $\E[f(T)]$.  However, a symmetric
  argument applied with the roles of the players interchanged implies
  that the value is no more than $\E[f(T)]$.
  
  Suppose that $M$ is an optimal strategy for the first player.  We
  have seen that when both players use $M$, the expected payoff is
  $\E[f(T)]=E(\emptyset,\emptyset)$.  Since $M$ is optimal for
  player~I, it follows that when both players use $M$ player~II always
  plays optimally (otherwise, player~I would gain an advantage, since
  she is playing optimally).  This means that $M(\emptyset,\emptyset)$
  is an optimal first move for player~II, and therefore every optimal
  first move for player~I is an optimal first move for player~II.  Now
  note that the game started at any position is equivalent to a
  selection game.  We conclude that every optimal move for one of the
  players is an optimal move for the other, which completes the proof.
\end{proof}

If $f$ is identically zero, then all strategies are optimal.  However,
if $f$ is \textit{generic\/} (meaning that all of the values $f(S_1)$
for different subsets $S_1$ of $S$ are linearly independent over
$\mathbb Q$), then the preceding argument shows that the optimal
choice of $s$ is always unique and that it is the same for both
players.  We thus have the following result:

\begin{prop}\label{p.generic}
  If $f$ is generic, then there is a unique optimal strategy and it is
  the same strategy for both players.  Moreover, when both players
  play optimally, the final $S_1$ is equally likely to be any one of
  the $2^n$ subsets of $S$.
\end{prop}

Theorem~\ref{t.rts} and Proposition~\ref{p.generic}
are in some ways quite surprising.  In the
baseball team selection, for example, one has to think very hard in
order to play the game optimally, knowing that at each stage there
is exactly one correct choice and that the adversary can capitalize
on any miscalculation. Yet, despite all of that mental effort by the
team captains, the final teams look no different than they would
look if at each step both captains chose players uniformly at
random.

Also, for purposes of illustration, suppose that there are only two players who
know how to pitch and that a team without a pitcher always loses.
In the alternating turn game, a captain can always wait to select a
pitcher until just after the other captain selects a pitcher.  In
the random-turn game, the captains must try to select the pitchers
in the opening moves, and there is an even chance the pitchers will
end up on the same team.

Theorem~\ref{t.rts} and Proposition~\ref{p.generic} generalize to
random-turn selection games in which the player to get the next turn
is chosen using a biased coin.  If player~I gets each turn with
probability $p$, independently, then the value of the game is $\E[f(T)]$, where
$T$ is a random subset of $S$ for which each element of $S$ is
in $T$ with probability $p$, independently.
For the corresponding statement of the proposition to hold, the notion
of ``generic'' needs to be modified. For example, it suffices to
assume that the values of $f$ are linearly independent over $\mathbb Q[p]$.
The proofs are essentially the same.  We leave it as an exercise to the
reader to generalize the proofs further so as to include the following
two games:

\spar{Restaurant Selection}
Two parties (with opposite food preferences) want to select a dinner
location.  They begin with a map containing $2^n$ distinct points
in $\R^2$, indicating restaurant locations.  At each step,
the player who wins a coin toss may draw a straight line that divides
the set of remaining restaurants exactly in half and eliminate all the
restaurants on one side of that line.  Play continues until one
restaurant $z$ remains, at which time player~I receives payoff
$f(z)$ and player~II receives $-f(z)$.

\spar{Balanced Team Captains}
Suppose that the captains wish
to have the final teams equal in size (i.e., there are $2 n$ players
and we want a guarantee that each team will have exactly $n$ players
in the end).  Then instead of tossing coins, the captains may shuffle a
deck of $2n$ cards (say, with $n$ red cards and $n$ black cards).  At
each step, a card is turned over and the captain whose color is shown
on the card gets to choose the next player.

\section{WIN-OR-LOSE SELECTION GAMES.}

We say that a game is a \textem{win-or-lose} game
if $f(T)$ takes on precisely two values, which we may as well
assume to be $-1$ and $1$.
If $S_1\subset S$ and $s\in S$, we say
that $s$ is \textem{pivotal} for $S_1$ if
$f(S_1\cup\{s\})\ne f(S_1\setminus\{s\})$.
A selection game is \textem{monotone} if $f$ is monotone; that is,
$f(S_1)\ge f(S_2)$ whenever $S_1\supset S_2$.
Hex is an example of a monotone, win-or-lose game. For such games,
the optimal moves have the following simple description:

\begin{lemma}\label{l.piv}
In a monotone, win-or-lose, random-turn selection game,
a first move $s$ is optimal if and only if $s$ is an element
of $S$ that is most likely to be pivotal for a random-uniform
subset $T$ of $S$.
When the position is $(S_1,S_2)$, the move $s$ in $S\setminus(S_1\cup S_2)$
is optimal if and only if $s$ is an element of $S\setminus(S_1\cup S_2)$
that is most likely to be pivotal for $S_1\cup T$, where
$T$ is a random-uniform subset of $S\setminus(S_1\cup S_2)$.
\end{lemma}

\noindent
The proof of the lemma is straightforward at this point and is
left to the reader.

For win-or-lose games, such as Hex, the players may stop making moves
after the winner has been determined, and it is interesting to
calculate how long a random-turn, win-or-lose, selection game will last
when both players play optimally.  Suppose that the game is a
monotone game and that, when there is more than one optimal move, the
players break ties in the same way.  Then we may take the point of
view that the playing of the game is a (possibly randomized) decision
procedure for evaluating the payoff function $f$ when the items are
randomly allocated.  Let $\vec x$ denote the allocation of the items,
where $x_i=\pm 1$ according to whether the $i$\th item goes to the
first or second player.  We may think of the $x_i$ as input bits,
and the playing of the game is one way to compute $f(\vec x)$.  The
number of turns played is the number of bits of $\vec x$ examined
before $f(\vec x)$ is computed.  We can use certain inequalities from the
theory of Boolean functions to bound the average length of play.

Let $I_i(f)$ denote the influence of the $i$\th bit on $f$ (i.e., the
probability that flipping $x_i$ will change the value of $f(\vec x)$).
The following inequality is from O'Donnell and Servedio
\cite{odonnell-servedio}:
\begin{multline}
 \sum_i I_i(f) =
 \E\left[\sum_i f(\vec x) x_i\right] =
 \E\left[f(\vec x) \sum_i x_i 1_{\text{$x_i$ examined}}\right]
 \le \text{(by Cauchy-Schwarz)}\\
 \sqrt{\E[f(\vec x)^2]\,\E\left[\left(\sum_{\text{$i$: $x_i$ examined}} x_i\right)^2 \right]} 
=
\sqrt{ \E\left[\left(\sum_{\text{$i$: $x_i$ examined}} x_i\right)^2 \right]} 
= \sqrt{\E[ \text{\# bits examined}]}\,.\label{eq:os}
\end{multline}
The last equality is justified by noting that
$\E[x_i\,x_j\, 1_{\text{$x_i$ and $x_j$ both examined}}]=0$
when $i\ne j$, which holds since conditioned on $x_i$ being examined
before $x_j$, conditioned on the value of $x_i$, and conditioned
on $x_j$ being examined, the expected value of $x_j$ is zero.
By~\eqref{eq:os} we have $$\E[\text{\# turns}] \geq \left[\sum_i I_i(f)\right]^2.$$
We shall shortly apply this bound to the game of Random-Turn Hex.

We mention one other inequality from the theory of Boolean functions, this
one due to O'Donnell, Saks, Schramm, and Servedio
 \cite[Theorem 3.1]{odonnell-saks-schramm-servedio}:
\begin{equation}\label{eq:osss}
\Var[f] \leq \sum_i \Pr[\text{$x_i$ examined}] I_i(f).
\end{equation}
For the random-turn game this implies that
$$ \E[\text{\# turns}] \geq \frac{\Var[f]}{\max_i I_i(f)}.$$

\section{RANDOM-TURN HEX.}
\label{s.hex}

\subsection*{Odds of winning on large boards and under biased play.}
In the game of Hex, the propositions discussed earlier imply that the probability that
player~I wins is given by the probability that there is a left-right
crossing in independent Bernoulli percolation
on the sites (i.e., when the sites are independently and randomly colored black or white).  One perhaps surprising
consequence of the connection to Bernoulli percolation is that, if
player~I has a slight edge in the coin toss and wins the coin toss
with probability $1/2+\eps$, then for any $\eps>0$, any $\delta>0$,
and any $r>0$
there is a strategy for player~I
that wins with probability at least $1-\delta$ on
the $L\times rL$ board, provided that $L$ is sufficiently large.

\subsection*{Random-Turn Hex on ordinary-size boards.}

Reisch proved in 1981 that determining which player has a winning
strategy on an arbitrarily shaped Hex board (with some of the sites
already colored in) is PSPACE-complete \cite{reisch}.  The online
community has, however, made a good deal of progress (much of it
unpublished) on solving the problem for smaller boards.  Jing Yang
\cite{yang} has announced the solution of Hex (and provided associated
computer programs) on boards of size up to $9\times 9$.  Hex is
usually played on an $11 \times 11$ board, for
which the optimal strategy is not yet known.

What is the situation with Random-Turn Hex?  We do not know if the
correct move in Random-Turn Hex can be found in polynomial time.  On
the other hand, for any fixed $\eps$ a computer can sample $O(L^4 \eps^{-2}
\log (L^4/\eps))$ percolation configurations (filling
in the empty hexagons at random) to estimate which empty site is most
likely to be pivotal given the current board configuration.  Except
with probability $O(\eps/L^2)$, the computer will pick a site that is
within $O(\eps/L^2)$ of being optimal.  This simple randomized strategy
provably beats an optimal opponent $50 - \eps$ percent of the time.

\begin{figure}[phtb]
\centerline{\includegraphics[width=.48\textwidth]{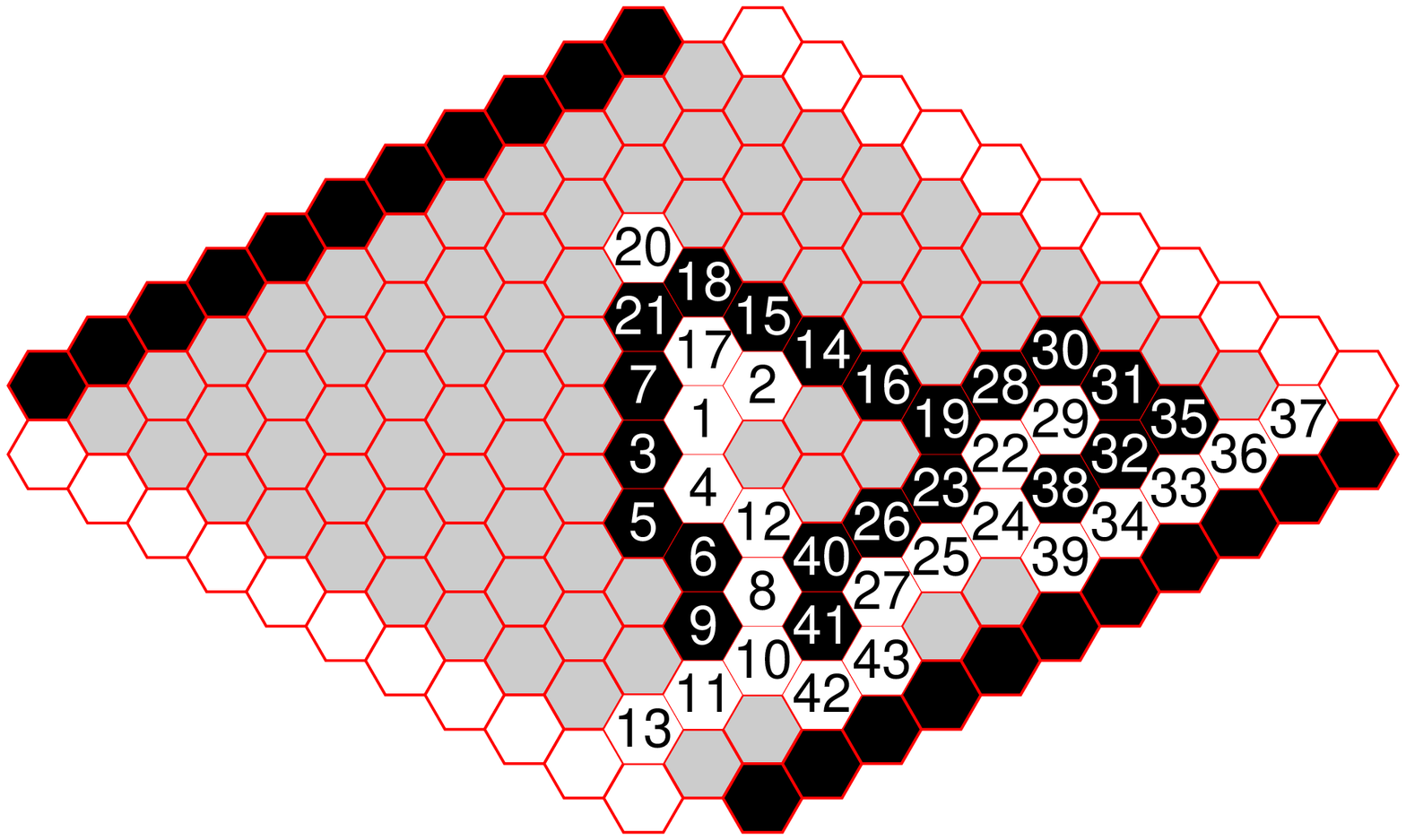}\hfill\includegraphics[width=.48\textwidth]{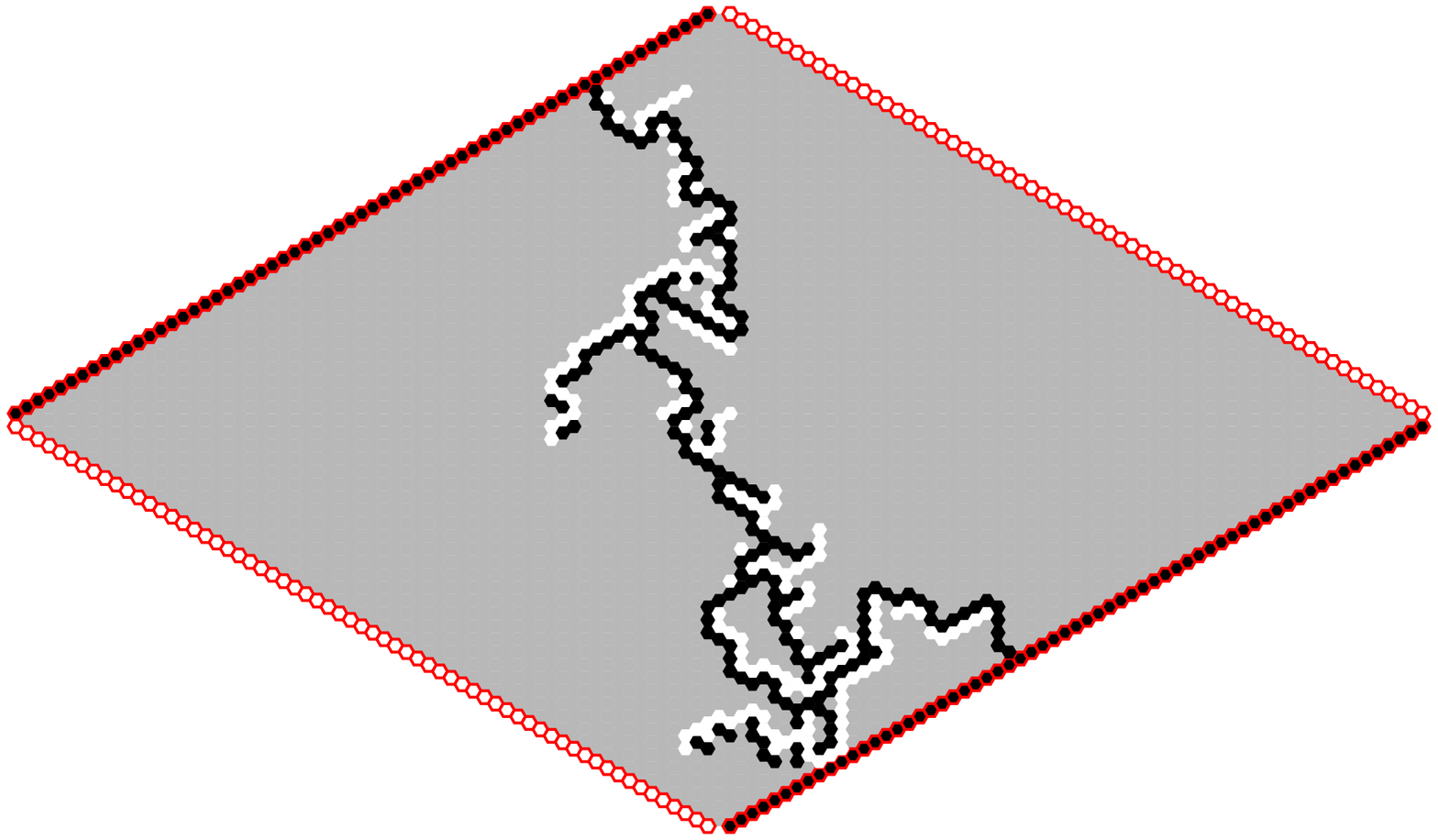}}
\caption{Random-Turn Hex on boards of size $11\times 11$ and $63\times 63$ under (near) optimal play.}
\label{fig:hex11}
\label{fig:hex63}
\end{figure}

\subsection*{Typical games under optimal play.}

What can we say about how long an average game of Random-Turn Hex will
last, assuming that both players play optimally?  (Here we assume
that the game is stopped once a winner is determined.)  If the side
length of the board is $L$, we wish to know how the expected length of
a game grows with $L$ (see Figure~\ref{fig:hex63} for games
on a large board).  Computer simulations on a variety of board sizes suggest
that the exponent is in the range $1.5$--$1.6$.  As far as rigorous bounds
go, a trivial upper bound is $O(L^2)$.  Since the game does not end
until a player has found a crossing, the length of the shortest
crossing in percolation is a lower bound, and empirically this
distance grows as $L^{1.1306\pm0.0003}$ \cite{MR1732786}, where the
exponent is known to be strictly larger than $1$.  We give a stronger
lower bound:
\begin{thm}\label{thm:hex-length}
  Random-Turn Hex under optimal play on an order $L$ board, when the
  two players break ties in the same manner, takes at least
  $L^{3/2+o(1)}$ time on average.
\end{thm}
\begin{proof}
To use the O'Donnell-Servedio bound~\eqref{eq:os}, we need to know the influence
that the sites have on whether or not there is a percolation crossing
(a path of black hexagons connecting the two opposite black sides).
The influence $I_i(f)$ is the probability that flipping site $i$ changes
whether there is a black crossing or a white crossing.  The ``4-arm
exponent'' for percolation is $5/4$ \cite{MR1879816}
(as predicted earlier in \cite{coniglio}), so 
$I_i(f)=L^{-5/4+o(1)}$ for sites $i$ ``away from the boundary,'' say in the
middle ninth of the region.  Thus $\sum_i I_i(f) \geq L^{3/4+o(1)}$,
so $\E[ \text{\# turns}] \geq L^{3/2+o(1)}$.
\end{proof}

\noindent
\textbf{Remark.} The function $\sum_i x_i 1_{\text{$x_i$ examined}}$
is the number of extra times that player~I wins the coin toss at the
time that the game terminates.  Naturally this is correlated with the
winner $f(\vec x)$ of the game, and for large board sizes this
correlation is noticeable.  Since the only inequality used in \eqref{eq:os}
was the Cauchy-Schwarz inequality, if we know $\sum_i I_i(f)$ and the
expected length of the game, we can determine how correlated the
winner is with who won most of the coin tosses.  Using more detailed
knowledge of how $I_i(f)$ behaves near the boundary, one can show that
for the standard lozenge-shaped Hex board, $\sum_i I_i(f)=L^{3/4+o(1)}$.
If the game lasts for $L^{>1.5}$ steps on average, this would imply
that the correlation between the winner of the game and the winner of
the majority of the coin tosses before the game is won tends to $0$ as
$L\to\infty$.

\vs

\noindent
\textbf{Remark.} The question of determining how many bits need to
be examined before one can decide whether or not there is a
percolation crossing arose in the context of dynamical percolation,
where the sites flip according to independent Poisson clocks
\cite{schramm-steif}.  Roughly speaking, if few bits need to be
examined, then it is easier for there to be exceptional times at which
there is an infinite percolating cluster in the plane
\cite{schramm-steif}.  One possible algorithm for determining whether
or not there is a crossing in an order $L$ region would be to follow
the black-white interfaces starting at the corners to see how they
connect, which exposes $L^{7/4+o(1)}$ hexagons in expectation
\cite{MR1879816}.  While this algorithm has the best currently
provable bound on the number of exposed hexagons, the ``play
Random-Turn Hex'' algorithm appears to do better (taking about
$L^{\text{$1.5$--$1.6$}}$ time).

\vs

\begin{figure}[phtb]
\centerline{\includegraphics[width=.35\textwidth]{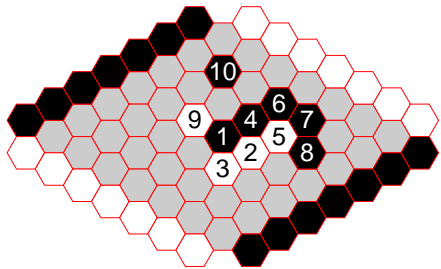}}
\caption{A rare occurrence---a game of Random-Turn Hex under (near) optimal play with a disconnected play.\label{fig:discon}}
\end{figure}

An optimally played game of Random-Turn Hex on a small board may
occasionally have a move that is disconnected from the other played
hexagons, as the game in Figure~\ref{fig:discon} shows.
But this is very much the exception rather than the rule.  For
moderate- to large-sized boards it appears that in almost every
optimally played game, the set of played hexagons remains a connected
set throughout the game (which is in sharp contrast to the usual game
of Hex).  We do not have an explanation for this phenomenon,
nor is it clear to us if it persists as the board size increases beyond
the reach of simulations.

\subsection*{Conformal (non)invariance of Random-Turn Hex.}

Combined with Smirnov's celebrated recent work in percolation
\cite{smirnov}, the connection between Random-Turn Hex and percolation
enables us to use Cardy's formula \cite{cardy} to approximate
player~I's probability of winning on very large boards of various
shapes, such as the $L\times r L$ rectangle.  The way this is done is
as follows.  Suppose that the players play the game on a very large
board that, when suitably scaled,
approximates some simply connected domain $D$ in the
complex plane, where the boundary of $D$ alternates colors at the
points $b_1,b_2,b_3,b_4$ with one black side between $b_1$ and $b_2$
and the other between $b_3$ and $b_4$.  By the Riemann Mapping
Theorem, there is an analytic function $\phi$ that bijectively maps
the domain $D$ to the upper half plane (i.e., a conformal bijection)
such that $\phi(b_1)=0$, $\phi(b_3)=1$, and $\phi(b_4)=\infty$.
Smirnov proved that in percolation on boards that are shaped like the
domain $D$, the probability of a black percolation crossing is a
function of $\phi(b_2)$ plus an error term that goes to zero as the
board size (the length of the shortest path connecting opposite sides)
goes to infinity.  This function of $\phi(b_2)$ is an
explicit formula (involving a hypergeometric function) that is known
as \textit{Cardy's formula}.  These crossing probabilities are known as a
``conformal invariant'' of percolation because they remain invariant
under conformal maps (up to error terms that go to zero as the board
size goes to infinity).  The connection between Random-Turn Hex and
percolation tells us that when the two players play optimally, the
probability of black winning is (up to an additive $o(1)$ error)
conformally invariant, and in particular is given by Cardy's formula.

\begin{figure}[phtb]
\centerline{\includegraphics[width=.58\textwidth]{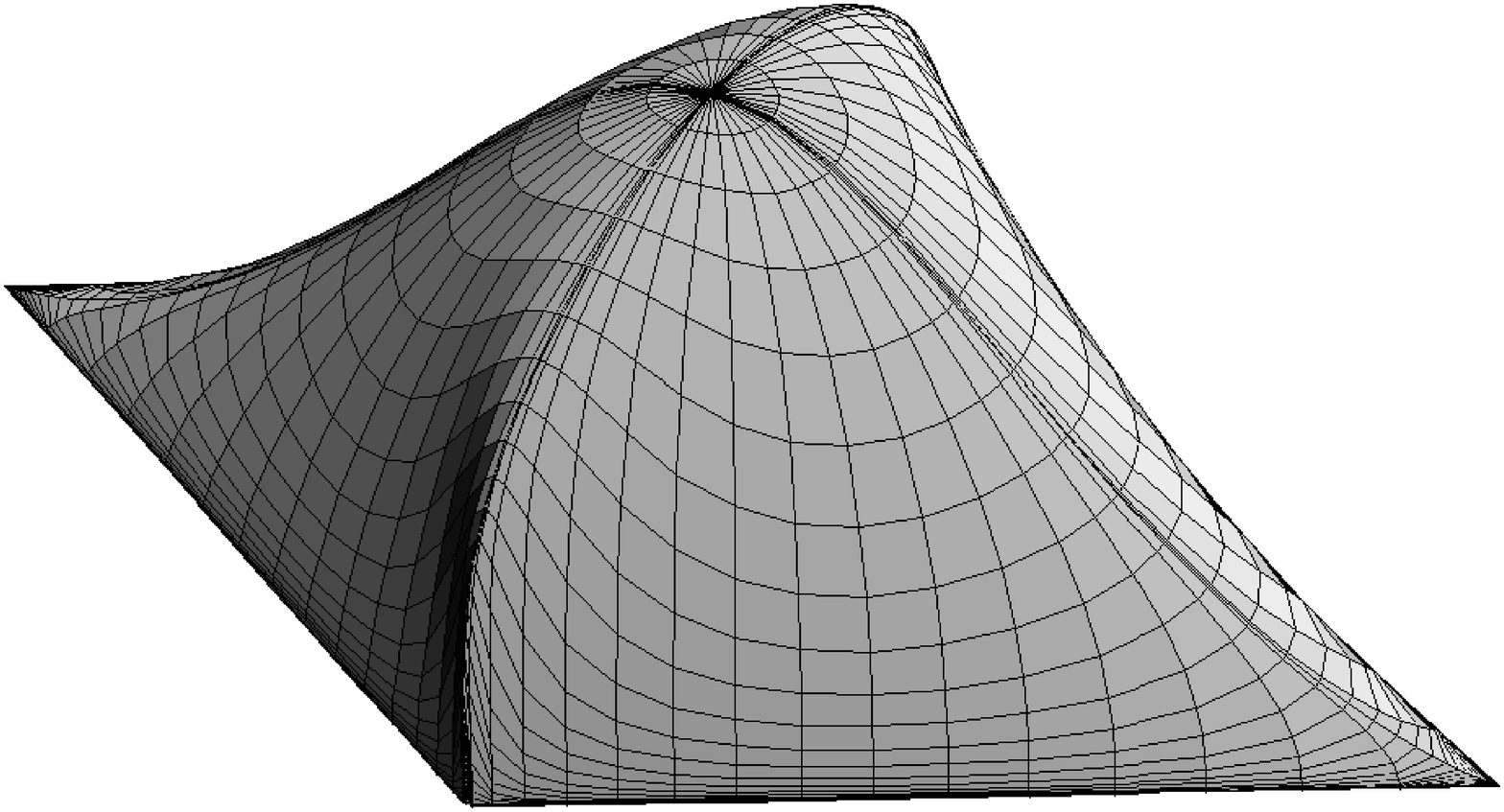}\hfill\includegraphics[width=.38\textwidth]{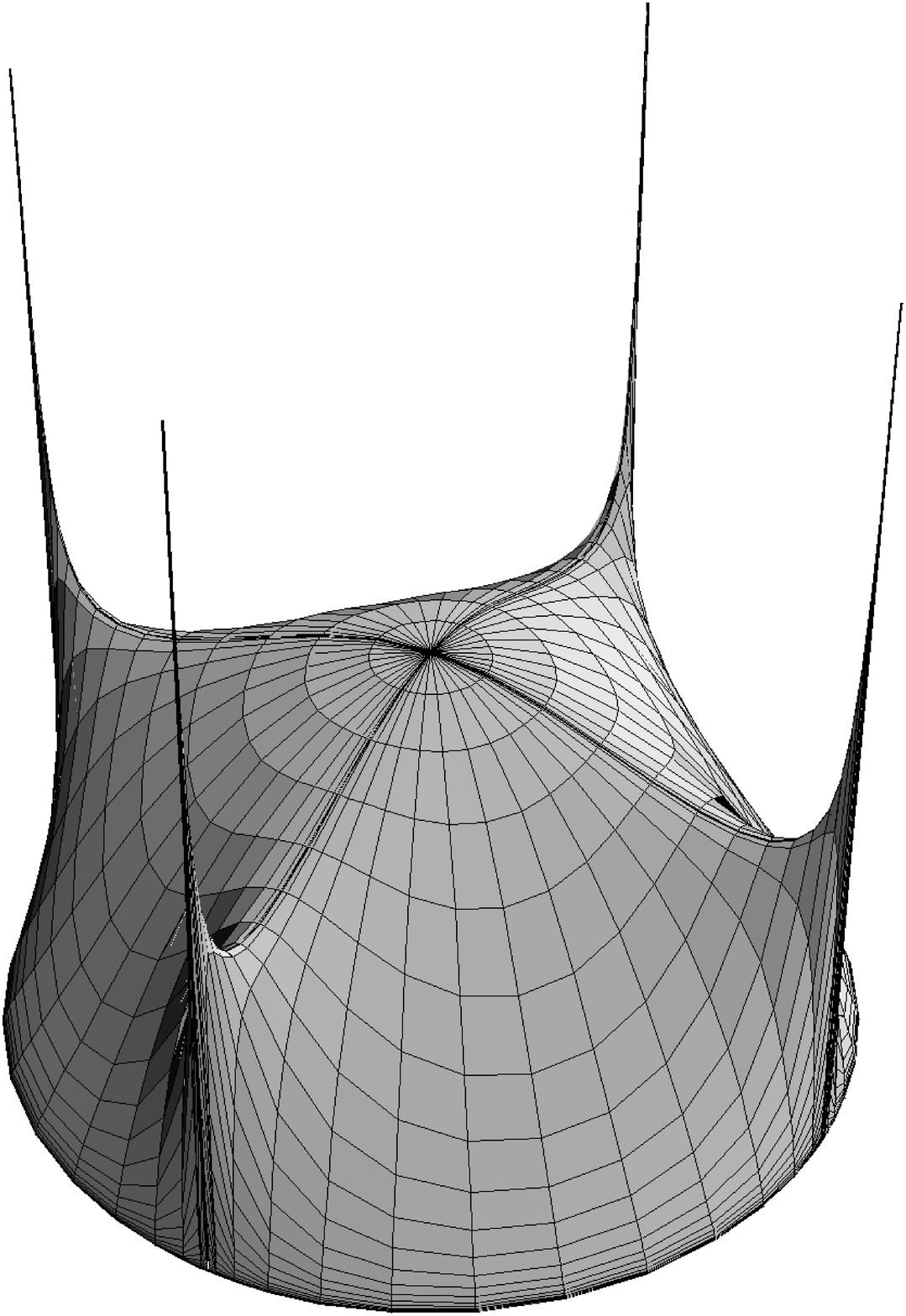}}
\caption{The location of the best first move in Random-Turn Hex is the site that is most likely to be pivotal for a percolation crossing.  Shown here are estimates of the probability that a site is critical when the board is the standard Hex board (on the left) and when the board is a disk (on the right), calculated using SLE \cite{schramm-wilson}.  For the standard board the best first move is near the center, while for the disk-shaped board the best first move is near one of the four points where the black boundary and white boundary meet.}
\label{fig:pivotal}
\end{figure}

However, the actual game play is \textit{not\/} conformally invariant,
and indeed, even the location of the first move is not conformally
invariant (see Figure~\ref{fig:pivotal}).  The reason for this is not
difficult to understand.  Recall that the first move is played at the
site that is most likely to be pivotal.  The probability that a site
is pivotal is not scale invariant (for larger boards, the probability
that a site is pivotal is smaller).  Since a general conformal map
scales different parts of the domain by different amounts, the
location of the site most likely to be pivotal is not conformally
invariant.

\section{THREE EXAMPLES ON TREES.}
\label{s.tree}

In this section, we study two different random-turn games based on trees.
Although these particular examples can be analyzed rather well, there are
other natural games based on trees that we do not know how to analyze.
We mention one such example at the end of the section (Recursive Three-Fold
Majority).

\subsection*{AND-OR trees.}

\begin{figure}[phtb]
\centerline{\includegraphics[width=\textwidth]{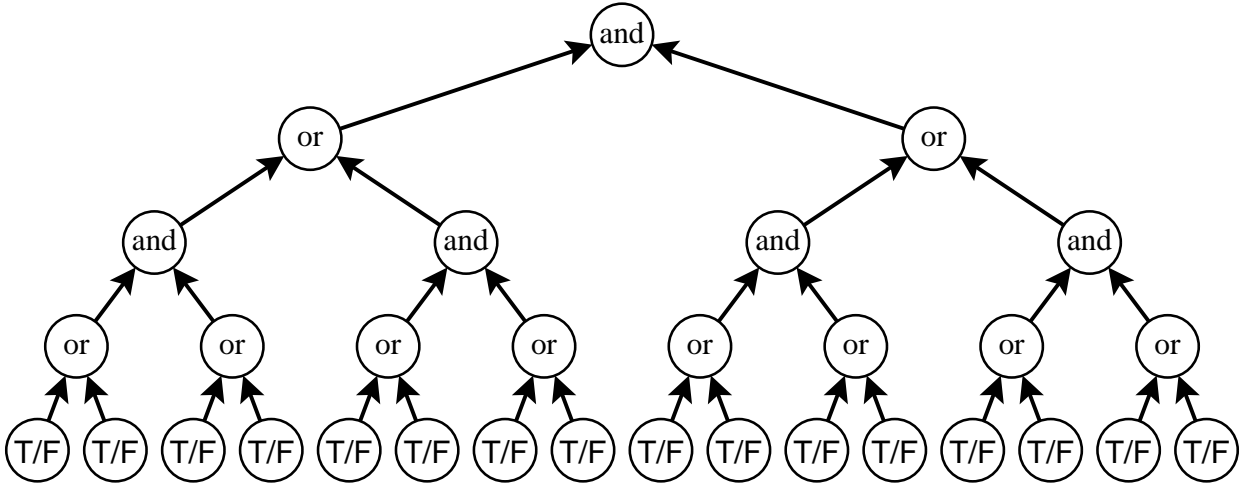}}
\caption{AND-OR tree of depth 4.  In Random-Turn AND-OR, the player who wins the coin toss sets the bit at a leaf node.  The value of the function determines the winner of the game.}
\label{fig:and-or}
\end{figure}

\newcommand{\Tr}{\textsf{T}\xspace}
\newcommand{\Fa}{\textsf{F}\xspace}
We believe that the expected length of game play in Random-Turn Hex
grows like a nonintegral power of the selection set size $|S|$ and that the
set of hexagons played has a random fractal structure.  The game of
Random-Turn AND-OR (see Figure~\ref{fig:and-or})
is a selection game for which we can actually prove
analogous statements.

The selection set $S$ is the set of leaves in a depth $h$ complete
binary tree.  We number the levels of the tree from $0$ (the root) to
$h$ (the leaves).  The leaves are treated as binary bits, with labels
of \Tr (``True'') or \Fa (``False'') if they are chosen by player~I or
II, respectively.  The internal nodes of the tree are also treated as
binary bits; the players cannot select the internal nodes, instead the
label of a bit on the $k$\th level of the tree ($k<h$) is the AND of
its two children if $k$ is even and the OR if $k$ is odd.  Player~I
wins if the root label is \Tr; player~II wins if it is \Fa.

When player~I always wins the coin toss with probability $p$ (which
need not be $1/2$), it is straightforward to estimate the probability
that she wins the game.  Let $q_k$ be the probability that a label at
the $k$\th level is \Tr.  Then $q_k = q_{k+1}^2$ (if $k$ is even) and
$q_k = 2 q_{k+1} - q_{k+1}^2$ (if $k$ is odd).  When $k$ is even,
$q_k = (2 q_{k+2} - q_{k+2}^2)^2$.  The map $q \mapsto (2q-q^2)^2$ has a
fixed point when $$q^4 - 4q^3 + 4q^2 -q = q(q-1)(q^2-3q+1) = 0$$
(i.e., if $q \in \{0,1,(3 \pm \sqrt{5})/2\}$); the fixed points in
$[0,1]$ are $0$, $1$, and $(3-\sqrt{5})/2 \approx .382$.  If $p$
is fixed and $h$ tends to $\infty$ along even integers, then the
probability player~I wins tends to $0$ if $p < (3-\sqrt 5)/2$
and to $1$ if $p > (3-\sqrt 5)/2$.

\begin{theorem}\label{th:andor}
  Consider an optimally played game of AND-OR on a level $h$ tree with
  coin-toss probability $p$.  If a move is played in some subtree $T$
  below some vertex $v$, then the succeeding moves are all played in
  $T$ until the label of $v$ is determined.  Moreover, the labels of
  the level $h-1$ vertices are determined in an order that is an
  optimally played game on the tree truncated at level $h-1$.
\end{theorem}
\begin{proof}
When a complete set of labels of the leaves (and hence all other
vertices) is fixed, we say a vertex on level $k$ is pivotal if it is a
pivotal element of the game on the tree truncated at level $k$.  A
vertex on level $k$ is pivotal if and only if its parent is pivotal
\textit{and\/} its sibling is \Tr (if $k$ is odd) or \Fa (if $k$ is
even).  Thus, a vertex is pivotal if and only if the path from that
vertex to the root has the property that every sibling of a vertex
along that path is \Tr or \Fa as the level is odd or even.  Note that
the event ``$v$ is pivotal'' is independent of labels on the subtree
below $v$.

We assume, inductively, that the statement of the theorem is true for
a tree with fewer than $h$ levels and for the first $t$ steps on a
level $h$ tree.  Before the $(t+1)$\th step is played, there are
either $0$ or $1$ partially determined vertices at level $h-1$ (that
is, vertices that are undetermined but for which a leaf below them has
been played).  If there are none, then the leaves that are most likely
to be pivotal are children of the level $h-1$ vertices that are most
likely to be pivotal, so any optimal $(t+1)$\th step also satisfies
the desired condition.  Suppose now that there is exactly one
partially determined vertex $v$ at level $h-1$, and call its labeled
child $x$ and its yet unlabeled child $y$.  Let $z$ be any unlabeled
leaf different from $y$ that might still affect the outcome of the
game.  Leaf $x$ must have been labeled at step $t$.  We need to show
that playing $y$ as the $(t+1)$\th move is better than playing $z$.

Let $r_t(u)$ be the conditional probability that a vertex $u$ is
pivotal given the labels that have been set prior to move $t$.  Let
$\bar p=p$ if $h$ is odd and $\bar p=1-p$ if $h$ is even ($\bar p$ is
the probability that the value of a leaf node does not determine its
parent).  Since $x$ was preferred over $z$, we have $r_t(x)\ge
r_t(z)$.  Now $r_t(z)$ is a nonnegative martingale, which is to say
that $\E[r_{t+1}(z)|\text{labels at time $t$}]=r_t(z)\geq 0$.  Since
$x$ got labeled the way it did with probability $\bar p$, we have
$\bar p r_{t+1}(z)\le r_t(z)\le r_t(x)=\bar p r_t(v)=\bar p
r_{t+1}(y)$.  Thus, we see that at step $t+1$, playing $y$ is at least
as good as playing $z$. We need only rule out the case of equality,
namely, $r_{t+1}(z)=r_{t+1}(y)$.

Let $a$ be the least common ancestor of $y$ and $z$, let $b_1$ be the
child of $a$ above $y$, and let $b_2$ be the child of $a$ above $z$.
Suppose that $x$ was the first move played in the subtree below $a$.
Regardless of which player got leaf $x$, the label of $a$ would not be
determined in turn $t$.  Thus, $z$ would still be pivotal with
positive probability, so the martingale argument would actually give
the strict inequality $\bar p r_{t+1}(z)< r_t(z)$ in place of the weak
inequality we used earlier.  Consequently, we get
$r_{t+1}(z)<r_{t+1}(y)$, as required.

On the other hand, suppose that $x$ was not the first move played
below $a$.  Then by induction, the $(t-1)$\th move was played below
$b_1$ and, again invoking the inductive hypothesis, we see that
$r_t(z)<r_t(x)$.  Thus, we get $r_{t+1}(z)<r_{t+1}(y)$ in this case as
well.
\end{proof}

Note that the Theorem~\ref{th:andor} completely characterizes the optimally
played games, whether or not the two players break ties in the same
way.  In fact, for every optimally played game there is an embedding
of the tree in the plane for which at each turn the played leaf is the
leftmost leaf that at that point has a positive probability to be
pivotal.

\begin{theorem}
  If $h$ is even and $p = (3-\sqrt 5)/2$, then the expected
  length of the game {\rm(}i.e., the expected number of labeled leaves{\rm)}
  is precisely $$\left(\frac{1+\sqrt 5}{2}\right)^h = (1.6180\ldots)^h.$$
\end{theorem}
\begin{proof}
  When all leaves are randomly labeled, a level $k$ vertex is $1$ with
  probability $p$ if $k$ is even and $1-p$ if $k$ is odd.  Under
  optimal play, once a vertex's label is determined, with
  probability $p$ it determines the label of its parent.  Thus,
  given that a vertex $v$ is labeled during an optimally played game,
  the expected number of labeled children of $v$ is $1+1-p$.
\end{proof}

Since we know the precise expected length of the game, we take a
moment to compare this with the lower bounds from the theory of
Boolean functions.  When the bits $x_i$ are true with probability $p$,
rather than taking $x_i=\pm1$ according to whether the bit is true or
false, it turns out to be more natural to take
$\Tr=\sqrt{(1-p)/p}$ and $\Fa=-\sqrt{p/(1-p)}$.  Then $\E[x_i]=0$
and $\E[x_i^2]=1$, so the O'Donnell-Servedio bound is still valid.
For AND-OR trees with $p=(3-\sqrt5)/2$, the influence of each
bit is $(\frac{-1+\sqrt5}{2})^h$, and there are $2^h$ bits, whence the
O'Donnell-Servedio bound gives
$$\E[\text{\# turns}] \geq
  (-1+\sqrt{5})^{2 h} = (1.5279\ldots)^h.
$$

The O'Donnell-Saks-Schramm-Servedio bound~\eqref{eq:osss}
generalizes to the
$p$-biased case~\cite[Theorem 3.1]{odonnell-saks-schramm-servedio}
and yields
$$\E[\text{\# turns}] \geq \frac{1}{(\frac{-1+\sqrt5}{2})^h} =
 \left(\frac{1+\sqrt 5}{2}\right)^h = (1.6180\ldots)^h.$$
This lower bound is exactly tight for AND-OR trees.

\subsection*{Shannon Switching Game.}

We now discuss a game for which we can prove that the set of moves is
always connected under optimal play.  In some ways, the analysis is
similar to the analysis of AND-OR games.

The \textit{Shannon Switching Game\/} is played by two players (named Cut and
Short) on a graph with two distinguished vertices.  When it is Cut's
turn, he may delete any unplayed edge of the graph, while Short may render an
edge immune to being cut.  Cut wins if he manages to disconnect the
two distinguished vertices, otherwise Short wins.  When the graph is
an $(L+1)\times L$ grid, with the vertices of the left side merged
into one distinguished vertex and the vertices on the other side
merged into another distinguished vertex, then the game is called
Bridg-It (also known as Gale after its inventor David Gale).  Oliver
Gross showed that the first player in Bridg-It has a (simple) winning
strategy, and then Lehman \cite{lehman} (see also \cite{mansfield})
showed how to solve the general Shannon Switching Game.  Here we
consider the random-turn version of the Shannon Switching Game.

Just as Random-Turn Hex is connected to site percolation on the
triangular lattice, where the vertices of the lattice (or, equivalently,
faces of the hexaognal lattice) are independently colored black or
white with probaiblity $1/2$, Random-Turn Bridg-It is connected to
bond percolation on the square lattice, where the edges of the square
lattice are independently colored black or white with probability
$1/2$.  We don't know the optimal strategy for Random-Turn Bridg-It,
but as with Random-Turn Hex, we can make a randomized algorithm that
plays near optimally.  Less is known about bond percolation than site
percolation, but it is believed that the crossing probabilities for
these two processes are asymptotically the same on
``nice'' domains~\cite{MR1230963}, %Langlands et al
so the probability that Cut wins in Random-Turn Bridg-It is
well approximated by the probability that a player wins in Random-Turn
Hex on a similarly shaped board.

Consider the Random-Turn Shannon Switching Game on a tree, where the
root is one distinguished vertex and the leaves are (collectively)
the other vertex.  Thus,
Short wins if at the end of the game there is a path from the
root to one of the leaves.  Under optimal play, the probability that Short wins
is just the probability that there is a path from the root to some
leaf when each edge is independently deleted with probability $1/2$.
For the complete ternary tree with $h$ levels, the probability that
Short wins under optimal play converges to $3-\sqrt{5}\approx0.764$
for large $h$ (see the proof of Theorem~\ref{thm:ternary}).
For the complete binary tree with $h$ levels Short
wins with probability $\sim 4/h$, so we define the ``enhanced binary
tree'' to be the tree where the root has $\lfloor h \log 2/2\rfloor$
children each of whom fathers a complete binary tree with $h-1$
levels.  Then the Random-Turn Shannon Switching Game on this enhanced
binary tree is an approximately fair game when $h$ is large, meaning
that Short wins with probability $1/2+o(1)$.  But what
are the optimal strategies, and how long does an average game last if
both players play optimally?  We will discuss these issues presently.

There is a natural coupling of a game with bond percolation on the tree.
The edges played by Short are the open edges, and those played by
Cut are the closed edges. The unplayed edges may be considered undecided.
We assume that the game terminates when the outcome is decided.

In the following result, a subtree below a vertex $v$ in a rooted tree
consists of $v$ and all the vertices and edges that are separated
from the root by $v$.

\begin{theorem}
  Consider the Random-Turn Shannon Switching Game on a tree of depth
  $h$ in which each internal node {\rm(}except possibly the root{\rm)} has $b$
  children and all the leaves are at level $h$.  Under optimal play,
  the set of moves played by Short forms a connected set of edges
  containing the root, and the set of moves played by Cut are all
  adjacent to the edges played by Short.  At any position, the set of
  optimal moves {\rm(}for either player{\rm)} consists of all the edges closest
  to the leaves among the undecided edges adjacent to shorted edges.
  Consequently, whenever an edge is played in some subtree $T$ below
  some vertex, the subsequent moves are played in $T$ until all the
  leaves of $T$ are disconnected from the root by cut edges or until
  Short wins.
\end{theorem}

\begin{proof}
After Short and Cut have made some moves, define the residual tree to
be the graph formed from the original tree by deleting the subtrees
disconnected from the root by cut edges and contracting
edges that have been shorted.
Consider some position of an optimally played game.
By induction, we assume that in the residual tree corresponding to
the current position
the direct descendants of the root are themselves roots of
regular $b$-ary subtrees of various depths whose leaves are
leaves of the original tree.
(The base of the induction is clear.)
Suppose that $e=[x,y]$ is an edge in the residual tree, where
$x$ is closer to the root than $y$.
The probability that $e$ is pivotal for bond percolation
from the root to the leaves (i.e., given the status of the
other edges, there is a connection from the root to the leaves if $e$
is uncut and there is no connection if $e$ is cut)
is precisely the probability that $y$ is connected to the leaves
(by uncut edges),
$x$ is connected to the root, and every open path from the root
to the leaves passes through $e$.
This is clearly maximized by edges in the residual tree adjacent to the root.
Among the edges adjacent to the root, this is maximized by
those closest to the leaves, because the probability
of having an open path to the root in Bernoulli percolation on
the regular $b$-ary tree is monotone decreasing in the depth of the tree.
By Lemma~\ref{l.piv}, these correspond to the optimal moves.
When one of these edges is shorted or cut, the stated
structure of the residual tree is preserved. Thus, the induction step is
established. The theorem follows.
\end{proof}

Next, we consider how long an average optimally played game lasts.  In
the following two theorems we use the notation (common in computer
science) $\Theta(g)$ to denote a quantity that is bounded between $c_1
g$ and $c_2 g$, where $c_1$ and $c_2$ are positive universal constants.
(The related notation $O(g)$ denotes an expression that is upper
bounded by $c_2 g$, but which may or may not be bounded below.)
\begin{theorem}
  Random-Turn Switching on the enhanced binary tree of depth $h$ lasts
  on average $\Theta(h^2)$ turns under optimal play.
\end{theorem}
\begin{proof}
  For the complete binary tree of depth $h$ the expected number of
  vertices in the percolation component of the root is $h+1$.
  Accordingly, the expected number
  of played edges of the enhanced binary
  tree is $O(h^2)$.  For the lower bound, note that for the complete
  binary tree the expected number of vertices in the percolation
  component of the root
  conditional on there being no leaf in this component is $(1+o(1)) h$.
  (This quantity is easily calculated inductively, for example.) For the
  enhanced binary tree, there is a good chance that $\Theta(h)$
  subtrees of the root are explored, so that on average at least
  $\Theta(h^2)$ moves are played.
\end{proof}
\begin{theorem}
  Random-Turn Switching on the complete ternary tree of order $h$
  lasts on average $\Theta(h)$ turns under optimal play.
\label{thm:ternary}
\end{theorem}
\begin{proof}
  For the complete ternary tree of depth $h$ let $q_h$ be the
  probability that Cut wins, let $\mu_h$ be the expected number of
  explored (played) edges conditional on Cut winning, and let $\nu_h$
  be the expected number of explored edges conditional on Short winning.
  We have $q_0=0$, $\mu_0=0$ (say), and $\nu_0=0$.
  When Cut wins, the top three edges are always explored.
  Conditional on Cut winning, these top edges are open with probability
  $q_{h-1}/(1+q_{h-1})$, independently. Consequently, the following recursions hold:
\begin{align*}
q_{h+1} &= \left(\frac{1+q_h}2\right)^3, \\
\mu_{h+1} &= 3 + 3\frac{q_h}{1+q_h} \mu_h\,, \\
% \intertext{since the top three edges are always explored, and conditional upon extinction, they are each open with probability $q_h/(1+q_h)$}
% \nu_{h+1} &= \frac{1}{1-q_h^3}\left[ (1-q_h^3) \mu_h + (1+q_h+q_h^2) + (q_h+q_h^2+q_h^3)\mu_h\right] \\
%           &= \mu_h + \frac{1}{1-q_h} + \frac{q_h}{1-q_h}\mu_h.\\
\nu_{h+1} &= \frac{1-q_h}{1-q_h^3}\left[ 1+\nu_h\right] + \frac{q_h(1-q_h)}{1-q_h^3}\left[2+\mu_h+\nu_h\right] + \frac{q_h^2(1-q_h)}{1-q_h^3}\left[3+2\mu_h+\nu_h\right] \\
&= \nu_h + \frac{1+q_h+q_h^2-3 q_h^3}{1-q_h^3} + \frac{q_h+q_h^2-2 q_h^3}{1-q_h^3} \mu_h\,.
\end{align*}
We then have $q_h\to \sqrt 5-2$ ($<1/2$), $\mu_h=\Theta(1)$, $\nu_h=\Theta(h)$, and $\E[\text{\# turns}] = \Theta(h)$.
% \begin{align*}
% \lim_{h\to\infty} q_h &=\sqrt{5}-2 =: q\\
% \lim_{h\to\infty} \mu_h &= 3\frac{1+q}{1-2q} = 3+\frac{9}{\sqrt{5}} =: \mu\\
% \lim_{h\to\infty} \frac{\nu_h}h &= \frac{1+q \mu}{1-q} -3 q^3\frac{1+\mu}{1-q^3}\\
% \lim_{h\to\infty} \frac{\E[\text{\# turns}]}h &= (1-q)\lim_{h\to\infty} \frac{\nu_h}h = \frac{70-54/\sqrt{5}}{19} = \Theta(1). \qedhere
% \end{align*}
\end{proof}

\subsection*{Recursive Majority.}

Let $S_h=S$ be a subset of the leaves of the complete ternary tree of
depth $h$.  Inductively, let $S_j$ be the set of nodes at level $j$
such that the majority of the nodes at level $j+1$ under them is in
$S_{j+1}$.  The payoff function $f(S)$ for \textit{Recursive
Three-Fold Majority\/} is $-1$ if $S_0=\emptyset$ and $+1$ if
$S_0=\{\text{root}\}$.  It seems that this random-turn game cannot be
analyzed using the methods we have used for the Random-Turn AND-OR
game or the Random-Turn Shannon Switching Game.  In particular, we do not know
how long the game takes (on average) when played optimally.

\section{OPEN PROBLEMS.}

We recall here some of the open problems raised earlier in the
article.  It would be interesting to know the expected length of a
game of optimally played Random-Turn Hex and whether or not the true
optimal move (not just near-optimal) can be found efficiently.  The
algorithm ``play Random-Turn Hex'' appears to be an efficient
algorithm (in terms of expected number of input bits examined) for
determining whether or not there is a percolation crossing, but is there
an asymptotically more efficient algorithm?  It would also be interesting
to know, in optimally played Random-Turn Hex on large boards, whether
or not the fraction of disconnected moves is asymptotically zero.
Finally, we do not know how long Random-Turn Ternary Recursive
Majority takes, or whether or not it is the most efficient algorithm
(in terms of expected number of input bits read) for evaluating
Recursive Ternary Majority.

\subsection*{ACKNOWLEDGMENTS.}

We are grateful to Wendelin Werner for suggesting to us that there
should be a version of Hex manifesting some conformal invariance in
the scaling limit.  This was our original motivation for considering
the random-turn version of the game.

Two of us (Peres and Sheffield) were supported in part by NSF grants \#DMS-0244479 and \#DMS-0104073.

\hrule
\vs \noindent
\textbf{YUVAL PERES} works in probability theory and fractal geometry,
especially on problems involving Brownian motion, probability on
trees, percolation, and projections of fractals.  After graduating
from the Hebrew University, Jerusalem, in 1990, he has taught at
Stanford, Yale, and the Hebrew University.  He now teaches in the
Statistics and Mathematics Departments at the University of
California, Berkeley.

\noindent\textit{Department of Statistics; 367 Evans Hall; University of California; Berkeley, CA 94720}%-3860}

\noindent{\texttt{http://www.stat.berkeley.edu/\char126peres/}}

\vs \noindent
\textbf{ODED SCHRAMM} was born in December 1961 in Jerusalem, Israel.  He studied
at the Hebrew University (BSc and MSc) and completed his Ph.D. at
Princeton University under the direction of William Thurston.  After a
two-year postdoc at U.C. San Diego, he joined the Weizmann Institute of
Science and remained there until 1999, at which time he moved to
Redmond, Washington, to join the Theory Group of Microsoft Research.  Among
the prizes he has received are the Salem Prize (2001), the Clay Research Award
(2002), the Henri Poincar\'e Prize (2003) and the Loeve Prize (2003).

\noindent\textit{Microsoft Research; One Microsoft Way; Redmond, WA 98052}

\noindent{\texttt{http://research.microsoft.com/\char126schramm/}}

\vs \noindent
\textbf{SCOTT SHEFFIELD} earned his Ph.D. from Stanford in 2003 and his B.A.
and M.A. degrees from Harvard in 1998.  He completed postdoctoral
research at Microsoft Research and U.C. Berkeley and is currently an
assistant professor in the Courant Institute at New York University,
where he teaches probability and mathematical finance.  He has held
numerous internships and short-term research positions in government
and industry, and his primary research interests are mathematical
physics and probability theory.

\noindent\textit{Courant Institute; 251 Mercer Street; New York, NY 10012}

\noindent{\texttt{http://www.cims.nyu.edu/\char126sheff/}}

\vs\noindent \textbf{DAVID B. WILSON} works in probability theory,
especially in statistical physics and random walks.  He studied at MIT
where he earned three S.B. degrees in 1991 and his Ph.D. in 1996 under
the direction of Jim Propp.  After postdocs at U.C. Berkeley, DIMACS
at Rutgers, and the Institute for Advanced Study, he joined the Theory
Group of Microsoft Research.

\noindent\textit{Microsoft Research; One Microsoft Way; Redmond, WA 98052}

\noindent{\texttt{http://dbwilson.com}}

\end{document}